\documentclass[12pt,oneside]{article}
\usepackage[a4paper, top=3cm, bottom=3cm, left=2.5cm, right=2.5cm]{geometry}
\usepackage{amsmath,amsfonts,amsthm,amssymb, graphicx,color}
\newtheorem{theorem}{Theorem}[section]

\newtheorem{lemma}{Lemma}[section]
\newtheorem{corollary}{Corollary}[section]

\newtheorem{remark}{Remark}[section]

\numberwithin{equation}{section}
\allowdisplaybreaks
\begin{document}

\title{ A local pointwise inequality for a biharmonic equation with negative exponents
\thanks{This work was partly supported by National Natural Science Foundation of China(NSFC)(No. 11971253), NSFC of Fujian(No. 2021J011101) and Fund of Key Laboratory of Financial Mathematics of Fujian Province University (Putian University).} }
\author{{\small Fan Chen$^{1,2}$, Jianqing Chen$^{1}$, Qihua Ruan$^{2}$\thanks{%
Corresponding author.  E-mail address: ruanqihua@163.com }
  } \\
%EndAName
{\small 1.College of Mathematics and Statistics $\&$ FJKLMAA,}\\
{\small Fujian Normal University, Qishan Campus, Fuzhou, 350117, P. R. China}\\
{\small 2.School of Mathematics and Finance, Putian University,}\\
{\small Putian, 351100, P. R. China}
 {\noindent\small }
}
\date{}
\maketitle

\begin{abstract}
In this paper, we are inspired by Ng\^{o}, Nguyen and Phan's \cite{QA} study of the pointwise inequality for positive $C^{4}$-solutions of biharmonic equations with negative exponent by using the growth condition of solutions. They propose an open question of whether the growth condition is necessary to obtain the pointwise inequality. We give a positive answer to this open question. We establish the following local pointwise inequality
$$-\frac{\Delta u}{u}+\alpha\frac{|\nabla u|^{2}}{u^{2}}+\beta u^{-\frac{q+1}{2}}\leq\frac{C}{R^{2}}$$
for positive $C^{4}$-solutions of the biharmonic equations with negative exponent
$$-\Delta^{2}u=u^{-q} \ \rm in \ \emph B_{\emph R}$$
where $B_{R}$ denotes the ball centered at $x_{0}$ with radius $R$, $n\geq3$, $q>1$, and some constants $\alpha\geq0$, $\beta>0$, $C>0$.


\numberwithin{equation}{section}

\end{abstract}

\section{Introduction and main results}
\label{laber}

In this paper, we consider a pointwise inequality for positive solutions of biharmonic equations with negative exponents
\begin{eqnarray}\label{1.7}
-\Delta^{2}u=u^{-q} \ \rm in \ \emph B_{\emph R},
\end{eqnarray}
where $B_{R}$ denotes the ball centered at $x_{0}$ with radius $R$, $n\geq3$, $q>1$.

Let us briefly describe the geometric background of this equation. Let $g=(g_{ij})$ be the standard Euclidean metric on $\mathbb{R}^{n}$, with $g_{ij}=\delta_{ij}$, where $\delta_{ij}=1$ for $i=j$, $\delta_{ij}=0$ otherwise. Let $\overline{g}=u^{\frac{4}{n-4}}g$, $n\neq4$, be a second metric derived from $g$ by the positive conformal factor $u: \mathbb{R}^{n}\rightarrow\mathbb{R}$. Then $u$ satisfies
$$\Delta^{2}u=\frac{n-4}{2}Q_{\overline{g}}u^{\frac{n+4}{n-4}},$$
where $Q_{\overline{g}}$ is the scalar curvature of $\overline{g}$, see \cite{Lai,Choi}. We know that if $n=3$, and $Q_{\overline{g}}>0$ is constant, then the equation \eqref{1.7} can be obtained by scaling for $q=7$. Recently, Lai and Ye \cite{LY} showed that if $0<q\leq1$, then equation \eqref{1.7} admits no entire smooth solution. For any $q>1$, there exist radial solutions to \eqref{1.7}, see \cite{McKenna,Kusano}.

We know that the pointwise differential inequality has exerted a great influence on the theory of  elliptic partial differential equations. It can be used to solve some famous conjectures and  open problems, and also can be applied to the existence theory of solutions of nonlinear  partial differential equations. The following  pointwise differential inequality was proved by Modica \cite{Modica}, which is one of the main techniques to solve De Giorgi's \cite{Giorgi} conjecture for the Allen-Cahn equation.
\begin{theorem}{\rm(\cite{Modica})}\label{theorem 1.0}
Let $F\in C^{2}(\mathbb{R})$ be a nonnegative function and $u$ be a bounded entire solution of
\begin{eqnarray}\label{1.0}
\Delta u=F'(u)\ \rm in \ \mathbb{R}^{n}.
\end{eqnarray}
Then
\begin{eqnarray}\label{1.01}
|\nabla u|^{2}\leq2F(u).
\end{eqnarray}
\end{theorem}

When $F'(u)=u(u^{2}-1)$, the equation \eqref{1.0} is known as the Allen-Cahn equation. The so-called De Giorgi's conjecture is a monotonic classification of the entire solution of the Allen-Cahn equation in one direction. To be more precise, assume that $u\in C^{2}(\mathbb{R}^{n})$ be a solution to the Allen-Cahn equation satisfying $\partial_{x_{n}}u\geq0$. Then the level sets $\{u=\lambda\}$ must be hyperplanes, for $n\leq8$. We know that De Georgi's conjecture has advanced considerably, having been fully established in dimensions $n=2$ by Ghoussoub and Gui \cite{Ghoussoub} and for $n=3$ by Ambrosio and Cabr\'{e} \cite{Ambrosio}. A celebrated result by Savin \cite{Savin} established its validity for $4\leq n\leq 8$ under an extra assumption that
$$\lim \limits_{x_{n}\rightarrow \pm\infty}u(x',x_{n})=\pm1.$$
For $n\geq 9$, the monotone solutions constructed by  Pino, Kowalczyk and Wei \cite{PKW}, who provide some nontrivial stable solutions.

As far as we know, a similar inequality as \eqref{1.01} corresponding to the fourth-order equation  is not proved. But some similar results were proved by Aghajani, Cowan and R$\breve{a}$dulescu \cite{ACR} in a bounded domain. There are many other well known pointwise inequalities for fourth-order Lane-Emden equations with positive exponents, i.e.,
\begin{eqnarray}\label{1.1}
\Delta^{2}u=u^{q} \ \rm in \ \mathbb{R}^{n}.
\end{eqnarray}
We refer interested readers to \cite{Wei,Mostafa2,Souplet} and references therein.

There is a big difference between the fourth-order elliptic equations with positive exponents and negative exponents. It is well known that the negative exponent of the nonlinear term is more challenging than the positive exponent,  since there  exist solutions that grow linearly or superlinearly at infinity, see \cite{Ignacio,McKenna,Quoc,Choi} and references therein.

Recently, Guo and Wei \cite{Zongming} established the following pointwise inequality for the positive $C^{4}$-solution of equation \eqref{1.7},
\begin{eqnarray}\label{1.8}
\Delta u\geq\sqrt{\frac{2}{q-1}}u^{-\frac{q-1}{2}}.
\end{eqnarray}

Inspired by the above nice result, Ng\^{o}, Nguyen and Phan \cite{QA} proved the following pointwise differential inequality.
\begin{theorem}{\rm(\cite{QA})}\label{theorem 1.1}
Let $u>0$ be a $C^{4}$-solution to \eqref{1.7} in $\mathbb{R}^{n}$ with $n\geq3$. Assume that two positive constants $\alpha$, $\beta$ satisfy
\begin{eqnarray}\label{1.9}
\begin{cases}
\alpha\leq\frac{1}{2},\\
\beta\leq\sqrt{\frac{2}{q-1-4\alpha/n}},\\
q\geq3\alpha+\sqrt{9\alpha^{2}+(1-2\alpha)(1+16\alpha/n)}.
\end{cases}
\end{eqnarray}
Then the following pointwise inequality
\begin{eqnarray}\label{1.10}
\Delta u\geq\alpha\frac{|\nabla u|^{2}}{u}+\beta u^{-\frac{q-1}{2}}
\end{eqnarray}
holds in $\mathbb{R}^{n}$ for any solution $u$ satisfying the growth condition
\begin{eqnarray}\label{1.11}
u(x)=o(|x|^{\frac{2}{1-\gamma}})\ as \ |x|\rightarrow\infty,
\end{eqnarray}
where $\gamma$ is arbitrary in $[0,1)$ such that
\begin{eqnarray}\label{1.12}
\begin{cases}
\alpha+\frac{4\alpha(1-2\alpha)}{n}-3\alpha\gamma-\gamma^{2}+\gamma>0,\\
q-1-\frac{8\alpha}{n}-2\gamma>0.\\
\end{cases}
\end{eqnarray}
In particular, the pointwise inequality \eqref{1.10} always holds under the assumption \eqref{1.9} and
\begin{eqnarray}\label{1.13}
u(x)=O(|x|^{2}), as |x|\rightarrow\infty.
\end{eqnarray}
\end{theorem}

The growth conditions of \eqref{1.11} and \eqref{1.13} seem not too strict, since it is proved by Lai and Ye \cite{LY}  that the radial solutions of the equation \eqref{1.7} in $\mathbb{R}^{n}$ grow at most quadratically at infinity. However the inequality \eqref{1.8} was proved by Guo and Wei without any growth condition of the solutions. So Ng$\hat{o}$, Nguyen and Phan propose an open question of whether the growth condition is necessary to obtain the pointwise inequality \eqref{1.10} when $\alpha>0$.

In this paper,  We give a positive answer to this open question. Our method of the proof was mainly motivated by the above Ng$\hat{o}$, Nguyen and Phan's \cite{QA} result. Our main contribution is to find a new auxiliary function which is different from the one in \cite{QA}. When using the maximum principle to this auxiliary function, we do not need any growth condition of the solutions. Here is our main theorem:
\begin{theorem}\label{theorem 1.2}
Let $u>0$ be a $C^{4}$-solution to \eqref{1.7} on a ball $B_{R}$ centered at $x_{0}$ with radius $R$ in $\mathbb{R}^{n}(n\geq3)$. Assume that two constant $\alpha\geq0$ and $\beta>0$ satisfy
\begin{eqnarray}\label{1.14}
\begin{cases}
\alpha<\frac{-(n^{2}+4n-4)+\sqrt{(n^{2}+4n-4)^{2}+64n}}{16}(n\geq4)\ ,\ \alpha<\frac{-7+2\sqrt{13}}{4}(n=3), \\
\beta\leq\sqrt{\frac{2}{q-1-4\alpha/n}},\\
q\geq max\{\frac{8\alpha}{n}+3,3\alpha+\sqrt{9\alpha^{2}+(1-2\alpha)(1+16\alpha/n)}\}.
\end{cases}
\end{eqnarray}
Then
\begin{eqnarray}\label{1.15}
-\frac{\Delta u}{u}+\alpha\frac{|\nabla u|^{2}}{u^{2}}+\beta u^{-\frac{q+1}{2}}\leq\frac{C}{R^{2}}.
\end{eqnarray}
\end{theorem}

As a direct consequence of Theorem \ref{theorem 1.2}, by taking $\alpha=0$ and letting $R\rightarrow\infty$, we obtain the following corollary:
\begin{corollary}\label{Corollary 1.1}
Let $u>0$ be a $C^{4}$-solution to \eqref{1.7} in $\mathbb{R}^{n}$ with $n\geq3$. Then
\begin{eqnarray*}\label{1.15.1}
\Delta u\geq\sqrt{\frac{2}{q-1}}u^{-\frac{q-1}{2}}.
\end{eqnarray*}
\end{corollary}
\begin{remark}
This pointwise inequality was concluded by Guo and Wei in {\rm \cite{Zongming}}.
\end{remark}

As another direct consequence of Theorem \ref{theorem 1.2}, by letting $R\rightarrow\infty$, we obtain the following corollary:
\begin{corollary}\label{Corollary 1.2}
Under the assumptions of Theorem \ref{theorem 1.2}. Then
\begin{eqnarray*}\label{1.15.2}
\Delta u\geq\alpha\frac{|\nabla u|^{2}}{u}+\beta u^{-\frac{q-1}{2}}.
\end{eqnarray*}
\end{corollary}

\begin{remark}
This pointwise inequality is obtained when $\alpha>0$. Which answer the question in {\rm \cite{QA}}.
\end{remark}

\section{Proof of Theorem \ref{theorem 1.2}}
\label{laber}
For simplicity, setting $v=\rm{ln}\emph u$, $k=\frac{q+1}{2}$, we obtain the following Lemma:
\begin{lemma}\label{lemma 2.1}
Let $u$ be a smooth positive solution of the equation \eqref{1.7} on a ball $B_{R}$ centered at $x_{0}$ with radius $R$ in $\mathbb{R}^{n}$. Then
\begin{eqnarray}\label{2.1.1}
-\Delta^{2}v&=&2|\nabla v|^{2}\Delta v+4\nabla\Delta v\cdot\nabla v+(\Delta v)^{2}\\
&&+|\nabla v|^{4}+4\nabla^{2}v(\nabla v,\nabla v)+2\|\nabla^{2}v\|^{2}+e^{-2kv},\nonumber
\end{eqnarray}
where $\nabla^{2}v$ denotes the Hessian matrix of $v$ and $\|\nabla^{2}v\|$ denotes the Hilbert-Schmidt norm on matrices defined to be
\begin{eqnarray*}
\|\nabla^{2}v\|=(\sum \limits_{i,j} |v_{ij}|^{2})^{\frac{1}{2}}.
\end{eqnarray*}
\end{lemma}
$Proof \ of \  Lemma\  \ref{lemma 2.1}$.

Since $u=e^{v}$, we deduce that
\begin{eqnarray}\label{2.1.5}
\Delta^{2} u&=&e^{v}(2|\nabla v|^{2}\Delta v+2\nabla\Delta v\cdot\nabla v+(\Delta v)^{2}+\Delta^{2}v+|\nabla v|^{4})\\
&&+e^{v}(4\nabla^{2}v(\nabla v,\nabla v)+\Delta|\nabla v|^{2}).\nonumber
\end{eqnarray}
From Bochner formula
\begin{eqnarray}\label{2.1.6}
\Delta|\nabla v|^{2}=2\|\nabla^{2}v\|^{2}+2\nabla\Delta v\cdot\nabla v.
\end{eqnarray}
Substituting \eqref{1.7} and \eqref{2.1.6} into \eqref{2.1.5}, it follows that
\begin{eqnarray}\label{2.1.7}
-\Delta^{2}v&=&2|\nabla v|^{2}\Delta v+4\nabla\Delta v\cdot\nabla v+(\Delta v)^{2}\\
&&+|\nabla v|^{4}+4\nabla^{2}v(\nabla v,\nabla v)+2\|\nabla^{2}v\|^{2}+e^{-(q+1)v}.\nonumber
\end{eqnarray}
These complete the proof of this Lemma.\\

We define an auxiliary function $\omega$ by
\begin{eqnarray}\label{2.5}
\omega=-\Delta v+(\alpha-1)|\nabla v|^{2}+\beta e^{-kv},
\end{eqnarray}
where $\alpha\geq0$, $\beta>0$ are constants to be chosen later.

Now we divide the proof of Theorem \ref{theorem 1.2} into two cases.\\

Case I. When $n\geq4$. We obtain the following lemma:
\begin{lemma}\label{lemma 2.2}
Let $u$ be a smooth positive solution of the equation \eqref{1.7} on a ball $B_{R}$ centered at $x_{0}$ with radius $R$ in $\mathbb{R}^{n}(n\geq4)$. Assume that two constant $\alpha\geq0$ and $\beta>0$ satisfy \eqref{1.14}. Then
\begin{eqnarray}\label{2.2}
\Delta\omega&\geq&-2(\alpha+1)\nabla v\cdot\nabla\omega+(1+\frac{2\alpha}{n})\omega^{2}+\frac{2}{n}(2-4\alpha-n)\alpha\omega|\nabla v|^{2}\\
&&+\frac{2}{n}(1-2\alpha)\alpha|\nabla v|^{4}\nonumber.
\end{eqnarray}
\end{lemma}
$Proof \ of \  Lemma\  \ref{lemma 2.2}$.

We write the Laplacian of $\omega$
\begin{eqnarray}\label{2.6}
\Delta\omega=-\Delta^{2}v+(\alpha-1)\Delta|\nabla v|^{2}+\beta\Delta e^{-kv}.
\end{eqnarray}
Since
\begin{eqnarray}\label{2.9}
\Delta e^{-kv}=-ke^{-kv}\Delta v+k^{2}e^{-kv}|\nabla v|^{2}.
\end{eqnarray}
Combining \eqref{2.1.6}, \eqref{2.6}, \eqref{2.9} and Lemma \ref{lemma 2.1}, we obtain that
\begin{eqnarray}\label{2.11}
\Delta \omega&=&2|\nabla v|^{2}\Delta v+2(\alpha+1)\nabla\Delta v\cdot\nabla v+(\Delta v)^{2}+|\nabla v|^{4}\\
&&+4\nabla^{2}v(\nabla v,\nabla v)+2\alpha\|\nabla^{2}v\|^{2}+e^{-2kv}\nonumber\\
&&+k^{2}\beta e^{-kv}|\nabla v|^{2}-k\beta e^{-kv}\Delta v\nonumber.
\end{eqnarray}
From \eqref{2.5}, we rewrite the first term of the right side of \eqref{2.11} by
\begin{eqnarray}\label{2.12}
2|\nabla v|^{2}\Delta v=-2\omega|\nabla v|^{2}+2(\alpha-1)|\nabla v|^{4}+2\beta e^{-kv}|\nabla v|^{2},
\end{eqnarray}
the second term by
\begin{eqnarray}\label{2.13}
2(\alpha+1)\nabla\Delta v\cdot\nabla v&=&-2(\alpha+1)\nabla v\cdot\nabla \omega+4(\alpha^{2}-1)\nabla^{2}v(\nabla v,\nabla v)\\
&&-2(\alpha+1)k\beta e^{-kv}|\nabla v|^{2}\nonumber,
\end{eqnarray}
the third term by
\begin{eqnarray}\label{2.14}
(\Delta v)^{2}&=&\omega^{2}+(\alpha-1)^{2}|\nabla v|^{4}+\beta^{2} e^{-2kv}-2(\alpha-1)\omega|\nabla v|^{2}\\
&&-2\beta\omega e^{-kv}+2(\alpha-1)\beta e^{-kv}|\nabla v|^{2}\nonumber,
\end{eqnarray}
the last term by
\begin{eqnarray}\label{2.15}
-\beta ke^{-kv}\Delta v=\beta k\omega e^{-kv}-(\alpha-1)\beta ke^{-kv}|\nabla v|^{2}-\beta^{2}ke^{-2kv}.
\end{eqnarray}
Substituting \eqref{2.12}, \eqref{2.13}, \eqref{2.14} and \eqref{2.15} into \eqref{2.11}, we get that
\begin{eqnarray}\label{2.16}
\Delta\omega&=&-2(\alpha+1)\nabla v\cdot\nabla \omega+\omega^{2}-2\alpha\omega|\nabla v|^{2}+(k-2)\beta \omega e^{-kv}\\
&&+(k^{2}-(3\alpha+1)k+2\alpha)\beta e^{-kv}|\nabla v|^{2}+(1+\beta^{2}-\beta^{2}k)e^{-2kv}\nonumber\\
&&+2\alpha\|\nabla^{2}v\|^{2}+4\alpha^{2}\nabla^{2}v(\nabla v,\nabla v)+\alpha^{2}|\nabla v|^{4}\nonumber.
\end{eqnarray}
It follows that
\begin{eqnarray}\label{2.17}
\Delta\omega&=&-2(\alpha+1)\nabla v\cdot\nabla \omega+\omega^{2}-2\alpha\omega|\nabla v|^{2}+(k-2)\beta \omega e^{-kv}\\
&&+(k^{2}-(3\alpha+1)k+2\alpha)\beta e^{-kv}|\nabla v|^{2}+(1+\beta^{2}-\beta^{2}k)e^{-2kv}\nonumber\\
&&+2\alpha\|\nabla^{2}v+\alpha\nabla v\otimes\nabla v\|^{2}+\alpha^{2}(1-2\alpha)|\nabla v|^{4}\nonumber.
\end{eqnarray}
From the Cauchy-Schwarz inequality, we have that
\begin{eqnarray}\label{2.18}
2\alpha\|\nabla^{2}v+\alpha\nabla v\otimes\nabla v\|^{2}&\geq&\frac{2\alpha}{n}(\Delta v+\alpha|\nabla v|^{2})^{2}\\
&=&\frac{2\alpha}{n}\omega^{2}-\frac{4\alpha(2\alpha-1)}{n}\omega|\nabla v|^{2}-\frac{4\alpha}{n}\beta\omega e^{-kv}\nonumber\\
&&+\frac{2\alpha}{n}(2\alpha-1)^{2}|\nabla v|^{4}+\frac{4\alpha(2\alpha-1)}{n}\beta e^{-kv}|\nabla v|^{2}+\frac{2\alpha\beta^{2}}{n}e^{-2kv}.\nonumber
\end{eqnarray}
Combining \eqref{2.17} and \eqref{2.18}, we deduce that
\begin{eqnarray}\label{2.19}
\Delta\omega&\geq&-2(\alpha+1)\nabla v\cdot\nabla \omega+(1+\frac{2\alpha}{n})\omega^{2}+(\frac{4(1-2\alpha)}{n}-2)\alpha\omega|\nabla v|^{2}\\
&&+(k-2-\frac{4\alpha}{n})\beta \omega e^{-kv}+(\frac{2}{n}(1-2\alpha)^{2}-2\alpha^{2}+\alpha)\alpha|\nabla v|^{4}\nonumber\\
&&+(1+\beta^{2}-\beta^{2}k+\frac{2\alpha\beta^{2}}{n})e^{-2kv}\nonumber\\
&&+(k^{2}-(3\alpha+1)k+2\alpha+\frac{4\alpha}{n}(2\alpha-1))\beta e^{-kv}|\nabla v|^{2}\nonumber.
\end{eqnarray}
That is
\begin{eqnarray}\label{2.19.1}
\Delta\omega&\geq&-2(\alpha+1)\nabla v\cdot\nabla\omega+(1+\frac{2\alpha}{n})\omega^{2}+Q_{1}\alpha\omega|\nabla v|^{2}+Q_{2}\beta\omega e^{-kv} \\
&&+I_{1}\alpha|\nabla v|^{4}+I_{2}e^{-2kv}+I_{3}\beta|\nabla v|^{2} e^{-kv}\nonumber,
\end{eqnarray}
where
\begin{eqnarray}\label{2.20}
\begin{cases}
I_{1}=\frac{2}{n}(1-2\alpha)^{2}-2\alpha^{2}+\alpha,\\
I_{2}=1+\beta^{2}-\beta^{2}k+\frac{2\alpha\beta^{2}}{n}\\
\ \ \ \ =1+\frac{2}{n}\alpha\beta^{2}-\frac{q-1}{2}\beta^{2},\\
I_{3}=k^{2}-(3\alpha+1)k+2\alpha+\frac{4\alpha}{n}(2\alpha-1)\\
\ \ \ \ =\frac{q-1}{2}(\frac{q+1}{2}-\alpha)-\alpha(q-\frac{8\alpha}{n}+\frac{4}{n}),\\
Q_{1}=\frac{4(1-2\alpha)}{n}-2,\\
Q_{2}=k-2-\frac{4\alpha}{n}=\frac{q-3}{2}-\frac{4\alpha}{n}.
\end{cases}
\end{eqnarray}
It is noted that the second inequality of \eqref{1.14} guarantees $I_{2}\geq0$ and the third of \eqref{1.14} guarantees $I_{3},Q_{2}\geq0$.
\\Let $0\leq\alpha<\frac{1}{2}$, we have that
\begin{eqnarray}\label{2.42}
I_{1}&=&\frac{2}{n}(1-2\alpha)^{2}-2\alpha^{2}+\alpha\\
&=&(\frac{8}{n}-2)\alpha^{2}+(1-\frac{8}{n})\alpha+\frac{2}{n}\nonumber\\
&\geq&\frac{4-n}{n}\alpha+(1-\frac{8}{n})\alpha+\frac{2}{n}\nonumber\\
&=&-\frac{4}{n}\alpha+\frac{2}{n}\nonumber.
\end{eqnarray}
Therefore combined with \eqref{2.19.1} and \eqref{2.42}, these complete the proof of Lemma \ref{lemma 2.2}.\\
$Proof \ of \  Theorem\ref{theorem 1.2}\ in \ the \ Case$ I.

We choose a $C^{2}$ cut-off function $0\leq\eta=\eta(t)\leq1$ on $[0,+\infty)$, which is defined as follows:
\begin{eqnarray}\label{2.21}
\eta(t)=\begin{cases}
1,\ \ \ \ \ t\in[0,1],\\
>0,\ \ t\in(1,2),\\
0,\ \ \ \ \ t\in[2,+\infty).
\end{cases}
\end{eqnarray}
such that it satisfies that for some positive constant $C$,
\begin{eqnarray}\label{2.22}
\frac{|\nabla \eta|^{2}}{\eta}\leq C, \ |\Delta \eta|\leq C.
\end{eqnarray}
Let $\rho(x)$ be the distance function from $x_{0}$, for $R>1$, we define
 \begin{eqnarray}\label{2.23}
\varphi=\eta(\frac{\rho(x)}{R}),
\end{eqnarray}
and
 \begin{eqnarray}\label{2.24}
\omega_{R}=\varphi\omega.
\end{eqnarray}
Then we have
 \begin{eqnarray}\label{2.25}
\frac{|\nabla\varphi|^{2}}{\varphi}=\frac{|\nabla\eta|^{2}}{R^{2}\eta}\leq\frac{C}{R^{2}},\ |\Delta \varphi|\leq \frac{C}{R^{2}}.
\end{eqnarray}
If $\omega(x)\leq0$ for $\forall x\in \mathbb{R}^{n}$, then it is nothing to prove. Suppose that
$$M=\sup \limits_{B_{R}}\omega>0.$$
Noting that $\omega_{R}=0$ if $\rho(x)>2R$, then there exists $x_{R}\in B_{2R}$ such that
\begin{eqnarray}\label{2.26}
M_{R}=\max \limits_{B_{2R}}\omega_{R}=\omega_{R}(x_{R}).
\end{eqnarray}
Since $M>0$, then $M_{R}>0$ for $R>1$. According to the necessary conditions of the local maximum, we see that
\begin{eqnarray}\label{2.27}
\begin{cases}
\nabla \omega_{R}=0,\\
\Delta \omega_{R}\leq0,\ at\  x=x_{R}
\end{cases}
\end{eqnarray}
Which implies that
\begin{eqnarray}\label{2.28}
\nabla\omega=-\frac{\nabla\varphi}{\varphi}\omega,
\end{eqnarray}
and
\begin{eqnarray}\label{2.29}
0\geq\varphi\Delta\varphi+2\nabla\varphi\cdot\nabla\omega+\omega\Delta\varphi
=\varphi\Delta\omega-2\frac{|\nabla\varphi|^{2}}{\varphi}\omega+\omega\Delta\varphi,\nonumber
\end{eqnarray}
at $x=x_{R}$. Substituting \eqref{2.25} into \eqref{2.29}, we obtain that
\begin{eqnarray}\label{2.30}
\frac{C}{R^{2}}\omega\geq\varphi\Delta\omega,
\end{eqnarray}
at $x=x_{R}$. Using Lemma \ref{lemma 2.2}, we see that, at $x=x_{R}$
\begin{eqnarray}\label{2.31}
\frac{C}{R^{2}}\omega&\geq&-2(\alpha+1)\varphi\nabla v\cdot\nabla\omega+(1+\frac{2\alpha}{n})\varphi\omega^{2}+\frac{2}{n}(2-4\alpha-n)\alpha\varphi\omega|\nabla v|^{2}\\
&&+\frac{2}{n}(1-2\alpha)\alpha\varphi|\nabla v|^{4}\nonumber.
\end{eqnarray}
From \eqref{2.25}, Cauchy-Schwarz inequality and Young's inequality, we deduce that, at $x=x_{R}$
\begin{eqnarray}\label{2.32}
-2(\alpha+1)\varphi\nabla v\cdot\nabla\omega&=&2(\alpha+1)\omega\nabla\varphi\cdot\nabla v\\
&\geq&-2(\alpha+1)\omega|\nabla \varphi||\nabla v|\nonumber\\
&\geq&-\frac{1}{\varepsilon_{1}}(\alpha+1)^{2}\omega\frac{|\nabla\varphi|^{2}}{\varphi}-\varepsilon_{1}\varphi\omega|\nabla v|^{2}\nonumber\\
&\geq&-\frac{C}{R^{2}\varepsilon_{1}}(\alpha+1)^{2}\omega-\varepsilon_{1}\varphi\omega|\nabla v|^{2}.\nonumber
\end{eqnarray}
Assuming that $0\leq\alpha<\frac{1}{2}$. Substituting \eqref{2.32} into \eqref{2.31} and take $\varepsilon_{1}=\frac{4\alpha(1-2\alpha)}{n}$, we see that
\begin{eqnarray}\label{2.33.1}
\frac{C}{R^{2}}\omega&\geq&(1+\frac{2\alpha}{n})\varphi\omega^{2}-2\alpha\varphi\omega|\nabla v|^{2}+(\frac{2}{n}-\frac{4\alpha}{n})\alpha\varphi|\nabla v|^{4},
\end{eqnarray}
at $x=x_{R}$, where $C$ is a positive constant does not depend on $R$. Applying Young's inequality, we have that
\begin{eqnarray}\label{2.34}
-2\alpha\varphi\omega|\nabla v|^{2}\geq-\frac{\alpha}{\varepsilon_{2}}\varphi\omega^{2}-\varepsilon_{2}\alpha\varphi|\nabla v|^{4}.
\end{eqnarray}
Substituting \eqref{2.34} into \eqref{2.33.1} and take $\varepsilon_{2}=\frac{2}{n}(1-2\alpha)$, we obtain that, at $x=x_{R}$
\begin{eqnarray}\label{2.35}
\frac{C}{R^{2}}\omega\geq (1+\frac{2\alpha}{n}-\frac{\alpha n}{2(1-2\alpha)})\varphi\omega^{2}.
\end{eqnarray}
In order to guarantees $1+\frac{2\alpha}{n}>\frac{\alpha n}{2(1-2\alpha)}$, we need that
\begin{eqnarray}\label{2.46}
8\alpha^{2}+(n^{2}+4n-4)\alpha-2n<0.
\end{eqnarray}
We can obtain the two real roots of the equation $8\alpha^{2}+(n^{2}+4n-4)\alpha-2n=0$ with
\begin{eqnarray}\label{2.46.1}
\alpha_{1}=\frac{-(n^{2}+4n-4)-\sqrt{(n^{2}+4n-4)^{2}+64n}}{16}<0,
\end{eqnarray}
and
\begin{eqnarray}\label{2.46.2}
0<\alpha_{2}=\frac{-(n^{2}+4n-4)+\sqrt{(n^{2}+4n-4)^{2}+64n}}{16}<\frac{1}{2}.
\end{eqnarray}
Choosing $0\leq\alpha<\alpha_{2}$, from \eqref{2.35}, we obtain that
\begin{eqnarray}\label{2.38}
\omega=\varphi\omega\leq\max \limits_{B_{2R}}\varphi\omega\leq\frac{C}{R^{2}},\ \rm{in}\   \emph B_{\emph R} .
\end{eqnarray}
From \eqref{2.5}, which implies that
\begin{eqnarray}\label{2.38.1}
-\frac{\Delta u}{u}+\alpha\frac{|\nabla u|^{2}}{u^{2}}+\beta u^{-\frac{q+1}{2}}\leq\frac{C}{R^{2}},\ \rm in\ \emph B_{\emph R}.
\end{eqnarray}
These complete the proof of the Case I of Theorem \ref{theorem 1.2}.

Case II. When $n=3$. From \eqref{2.19.1}, the following lemma it is obvious to obtain.
\begin{lemma}\label{lemma 2.3}
Let $u$ be a smooth positive solution of the equation \eqref{1.7} on a ball $B_{R}$ centered at $x_{0}$ with radius $R$ in $\mathbb{R}^{3}$. Assume that two $\alpha\geq0$ and $\beta>0$ satisfy \eqref{1.14}. Then
\begin{eqnarray}\label{2.2.1}
\Delta\omega&\geq&-2(\alpha+1)\nabla v\cdot\nabla\omega+(1+\frac{2\alpha}{3})\omega^{2}-\frac{2}{3}(1+4\alpha)\alpha\omega|\nabla v|^{2}\\
&&+\frac{1}{3}(2\alpha^{2}-5\alpha+2)\alpha|\nabla v|^{4}\nonumber.
\end{eqnarray}
\end{lemma}
$Proof \ of \  Theorem\ref{theorem 1.2}\ in \ the \ Case$ II.

From \eqref{2.30} and Lemma \ref{lemma 2.3}, we see that
\begin{eqnarray}\label{2.2.2}
\frac{C}{R^{2}}\omega&\geq&-2(\alpha+1)\varphi\nabla v\cdot\nabla\omega+(1+\frac{2\alpha}{3})\varphi\omega^{2}-\frac{2}{3}(1+4\alpha)\alpha\varphi\omega|\nabla v|^{2}\\
&&+\frac{1}{3}(2\alpha^{2}-5\alpha+2)\alpha\varphi|\nabla v|^{4}\nonumber,
\end{eqnarray}
at $x=x_{R}$. From \eqref{2.28}, Cauchy-Schwarz inequality and Young's inequality, we obtain that
\begin{eqnarray}\label{3.1}
-2(\alpha+1)\varphi\nabla v\cdot\nabla\omega\geq-\frac{C}{R^{2}\varepsilon_{3}}(\alpha+1)^{2}\omega-\varepsilon_{3}\varphi\omega|\nabla v|^{2},
\end{eqnarray}
at $x=x_{R}$. Substituting \eqref{3.1} into \eqref{2.2.2} and take $\varepsilon_{3}=\frac{2}{3}(1+4\alpha)\alpha$, we have that, at $x=x_{R}$
\begin{eqnarray}\label{3.2}
\frac{C}{R^{2}}\omega\geq (1+\frac{2\alpha}{3})\varphi\omega^{2}-\frac{4}{3}(1+4\alpha)\alpha\varphi\omega|\nabla v|^{2}+\frac{1}{3}(2\alpha^{2}-5\alpha+2)\alpha\varphi|\nabla v|^{4}.
\end{eqnarray}\\
Applying Young's inequality, we see that, at $x=x_{R}$
\begin{eqnarray}\label{3.3}
-\frac{4}{3}(1+4\alpha)\alpha\varphi\omega|\nabla v|^{2}
\geq-\frac{4\alpha(1+4\alpha)^{2}}{9\varepsilon_{4}}\varphi\omega^{2}-\varepsilon_{4}\alpha\varphi|\nabla v|^{4}.
\end{eqnarray}
Assuming that $0\leq\alpha<\frac{1}{2}$, choose $\varepsilon_{4}=\frac{1}{3}(2\alpha^{2}-5\alpha+2)>0$, substituting \eqref{3.3} into \eqref{3.2}, we obtain that, at $x=x_{R}$
\begin{eqnarray}\label{3.4}
\frac{C}{R^{2}}\omega\geq (1+\frac{2\alpha}{3}-\frac{4\alpha(1+4\alpha)^{2}}{3(2\alpha^{2}-5\alpha+2)})\varphi\omega^{2}.
\end{eqnarray}
We notice that
\begin{eqnarray}
\frac{4\alpha(1+4\alpha)^{2}}{3(2\alpha^{2}-5\alpha+2)}=\frac{4\alpha(1+4\alpha)^{2}}{3(1-2\alpha)(2-\alpha)}<\frac{4\alpha(1+2)^{2}}{3(1-2\alpha)(2-\frac{1}{2})}=\frac{8\alpha}{1-2\alpha}.
\end{eqnarray}
In order to guarantees $1+\frac{2\alpha}{3}>\frac{4\alpha(1+4\alpha)^{2}}{3(2\alpha^{2}-5\alpha+2)}$, we need that
\begin{eqnarray}
1+\frac{2\alpha}{3}>\frac{8\alpha}{1-2\alpha}.
\end{eqnarray}
That is
\begin{eqnarray}\label{2.48}
4\alpha^{2}+28\alpha-3<0.
\end{eqnarray}
We can obtain the two real roots of the equation $4\alpha^{2}+28\alpha-3=0$ with
\begin{eqnarray}\label{2.49}
\alpha_{3}=\frac{-7-2\sqrt{13}}{4}<0,
\end{eqnarray}
and
\begin{eqnarray}\label{2.50}
0<\alpha_{4}=\frac{-7+2\sqrt{13}}{4}<\frac{1}{2}.
\end{eqnarray}
Choosing $0\leq\alpha<\alpha_{4}$, from \eqref{3.4}, we obtain that
\begin{eqnarray}\label{2.51}
\omega\leq\frac{C}{R^{2}},\ \rm{in}\   \emph B_{\emph R} ,
\end{eqnarray}
these complete the proof of the Case II of Theorem \ref{theorem 1.2}.

These complete the proof of Theorem \ref{theorem 1.2}.



%% References


\bibliographystyle{elsarticle-num}



%%%%%%%%%%%%%%%%%%%%%%%%%%%%%%%%%%%%%%%%%%%%%%%%%%%%%%%%%%%%%%%%%%%%%%%%%%%%%%%%%%%%%%%%%%%

\end{document}